# A matrix preconditioning framework for physics-informed neural networks based on adjoint method

Jiahao Song[1,2,3], Wenbo Cao[1,2,3], Weiwei Zhang[1,2,3,*]

[1] School of Aeronautics, Northwestern Polytechnical University, Xi'an 710072, China

[2] International Joint Institute of Artificial Intelligence on Fluid Mechanics, Northwestern Polytechnical University, Xi'an, 710072, China

[3] National Key Laboratory of Aircraft Configuration Design, Xi'an 710072, China

* Corresponding author. E-mail: aeroelastic@nwpu.edu.cn

**Abstract** Physics-informed neural networks (PINNs) have recently emerged as a popular approach for solving forward and inverse problems involving partial differential equations (PDEs). Compared to fully connected neural networks, PINNs based on convolutional neural networks offer advantages in the hard enforcement of boundary conditions and in reducing the computational cost of partial derivatives. However, the latter still struggles with slow convergence and even failure in some scenarios. In this study, we propose a matrix preconditioning method to improve the convergence of the latter. Specifically, we combine automatic differentiation with matrix coloring to compute the Jacobian matrix of the PDE system, which is used to construct the preconditioner via incomplete LU factorization. We subsequently use the preconditioner to scale the PDE residual in the loss function in order to reduce the condition number of the Jacobian matrix, which is key to improving the convergence of PINNs. To overcome the incompatibility between automatic differentiation and triangular solves in the preconditioning, we also design a framework based on the adjoint method to compute the gradients of the loss function with respect to the network parameters. By numerical experiments, we validate that the proposed method successfully and efficiently solves the multi-scale problem and the high Reynolds number problem, in both of which PINNs fail to obtain satisfactory results.

## 1 Introduction

In recent years, with advances in computational power and deep learning, deep neural networks have shown remarkable promise in solving forward and inverse problems involving partial differential equations (PDEs). Representative studies include physics-informed neural networks (PINNs) [1], the deep Galerkin method (DGM) [2], and the deep Ritz method (DRM) [3], among which PINNs have received the most attention. The core of PINNs is incorporating the PDE into the loss function to ensure that the neural network output satisfies physical constraints. In fact, this idea can be traced back to the 20th century, when Dissanayake and Phan-Thien [4] pioneered the use of neural networks to solve PDEs. Since it was proposed in 2018, PINNs have produced significant research outcomes in various fields, such as fluid mechanics [5-7], heat transfer [8, 9], materials [10, 11], and electromagnetic propagation [12].

PINNs typically approximate the solution of PDEs using a fully connected neural network and compute the partial derivatives in the PDEs by automatic differentiation [13]. Based on this framework, numerous methods have been proposed to improve the performance of PINNs. For example, Wang et al. [14] enhanced the high-frequency approximation ability of PINNs by embedding a Fourier feature layer after the input layer. Jagtap and Karniadakis [15] proposed a parallelizable domain decomposition method to improve the accuracy and efficiency of PINNs. Wu et al. [16] comprehensively investigated the effects of sampling methods on PINNs and developed two adaptive sampling strategies to reduce the number of collocation points. Cao and Zhang [17] established a connection between the ill-conditioning of PINNs and the Jacobian matrix of the PDE system. Wassing et al. [18] alleviated the difficulty of PINNs in shock capturing by introducing artificial viscosity into the Euler equations. Despite significant outcomes, PINNs still have limitations in certain aspects. First, in PINNs, conflicts between the PDE loss and the boundary condition loss may reduce the accuracy or even cause failure [19-21]. Although researchers have proposed several methods for the hard enforcement of boundary conditions [22-24] to avoid this issue, these methods still rely on empirical functions that affect the robustness of PINNs. Second, the cost of automatic differentiation increases exponentially with the order of partial derivatives, leading to the high computational cost for PINNs [25].

One way to address these two limitations is to approximate the solution of the PDE using a convolutional neural network (CNN). On the one hand, CNN allow hard constraints for boundary conditions to be enforced by padding [26, 27], similar to numerical methods, without introducing any additional parameters or functions. On the other hand, in CNN, the partial derivatives in PDEs are commonly computed via the finite difference, whose cost increases linearly with the order of partial derivatives, resulting in lower computational cost compared to automatic differentiation. Furthermore, discretizing the PDE using the finite difference allows for the integration of well-established techniques from traditional numerical methods, such as the upwind scheme and the TVD scheme, which may be crucial for improving the accuracy and convergence of the algorithm. According to the differences in mapping, existing studies on solving PDEs using CNN can be divided into two categories. The first is to establish a mapping from the solution at the current time step to the solution at the next time step. For example, Wandel et al. [28] developed a surrogate model for predicting three-dimensional incompressible flows. Ren et al. [27] combined CNN with long short-term memory (LSTM) to design a recurrent learning framework. Ranade et al. [29] proposed DiscretizationNet based on the finite volume discretization to predict steady-state problems. Their training mechanics are similar to the pseudo-time stepping. These studies are valuable for achieving stable long-term prediction. The second is to establish a mapping from space-time coordinates to the PDE solution, which is same as PINNs. Our study belongs to this category. Compared to the first category, there are fewer existing studies employing this mapping. Gao et al. [26] proposed the physics-informed geometry-adaptive convolutional neural network (PhyGeoNet) to solve PDEs defined on irregular domains. They used the coordinate transformation to convert the irregular physical domain to the regular reference domain for matching the input of CNN. Qi et al. [30] demonstrated the advantage of using CNN to approximate the PDE solution over the fully connected neural network by solving the 2D Shallow Water Equations. These studies have made significant contributions to enhancing PINNs. However, they do not address all the challenges faced by PINNs. For example, in some scenarios, the algorithm may exhibit slow convergence or even converge to incorrect solutions [17, 31]. Fast convergence is key to ensuring the reliability and generality of PDE solvers [32, 33]. Therefore, in this study, we focus on improving the convergence of existing methods.

We propose a matrix preconditioning method for PINNs. Matrix coloring and automatic differentiation are combined to compute the Jacobian matrix of the PDE system. Then, we construct the preconditioner using the incomplete LU factorization and utilize it to scale the PDE residual in the loss function of PINNs. To reduce memory and computational costs, we replace the computation of the inverse of the preconditioner with an iterative method. Meanwhile, the adjoint method is introduced to solve the issue that automatic differentiation is not supported in the iterative method involving sparse matrices, which results in the gradients of the loss function with respect to the network parameters being uncomputable. We validate the effectiveness of the proposed method by solving the multi-scale problem and the high Reynolds number problem, in both of which PINNs fail to obtain satisfactory results. The remainder of the paper is organized as follows. In Section 2, we first provide an overview of PINNs. Then, the matrix preconditioning method for PINNs is introduced in detail. In Section 3, we validate the effectiveness of the proposed method by numerical experiments. Finally, concluding remarks and direction for future research are presented in Section 4.

## 2 Methodology

2.1 Physics-informed neural networks

We consider the steady-state PDE defined on a domain $\Omega \subset \mathbb{R}^d$:

$$\begin{aligned}\mathcal{N}[u] &= 0, \quad x = (x_1,...,x_d) \in \Omega \\ \mathcal{B}[u] &= 0, \quad x \in \partial\Omega\end{aligned} \tag{1}$$

where $u$ denotes the solution of the PDE. $\mathcal{N}[\cdot]$ is the general differential operator which can include any combination of linear and nonlinear terms of spatial derivatives, and $\mathcal{B}[\cdot]$ is the boundary condition operator.

PINNs use the deep neural network (e.g., fully connected neural network, convolutional neural network) to approximate $u$. The network takes space coordinates as the input and outputs the approximate solution $u^{NN}(x;\theta)$. $\theta$ represent network parameters. When using CNN, the loss function of PINNs is:

$$\mathcal{L} = \frac{1}{N_r} \left\| \mathcal{N}[u^{NN}(x_r^i;\theta)] \right\|^2 \tag{2}$$

where $\{x_r^i\}_{i=1}^{N_r}$ are collocation points. Since the boundary conditions are enforced as hard constraints through padding [26, 27], Eq. (2) includes only the PDE loss. Conflicts [19-21] between multiple loss terms in the fully connected neural network are avoided. The partial derivatives in $\mathcal{N}[u^{NN}(x_r^i;\theta)]$ are computed via the finite difference, unless stated otherwise, we use the second-order central difference in this study. Figure 1 shows a schematic of PINNs, in which a CNN is employed to approximate $u$. We use the U-Net architecture [34], which connects shallow and deep features through skip connections, and also refer to the design by Hu et al. [35]. Specifically, the input feature map is first linearly transformed to 8 dimensions via 1×1 convolution kernels, while preserving its size. Then, the extended feature map is encoded and decoded sequentially. In the encoder, we use 4×4 convolutional kernels (stride = 2, padding = 1), where the feature map size is halved at each layer, and the number of channels doubles in turn. In the decoder, we use bicubic upsample with the scale factor of 2, followed by 3×3 convolutional kernels (stride = 1, padding = 1). The feature map size doubles at each layer, and the number of channels is halved in turn. Moreover, channel attention [36] is introduced to enhance the learning capability of the neural network. Finally, to obtain the output, the feature map is reduced in dimensionality via the linear transform while maintaining its size.

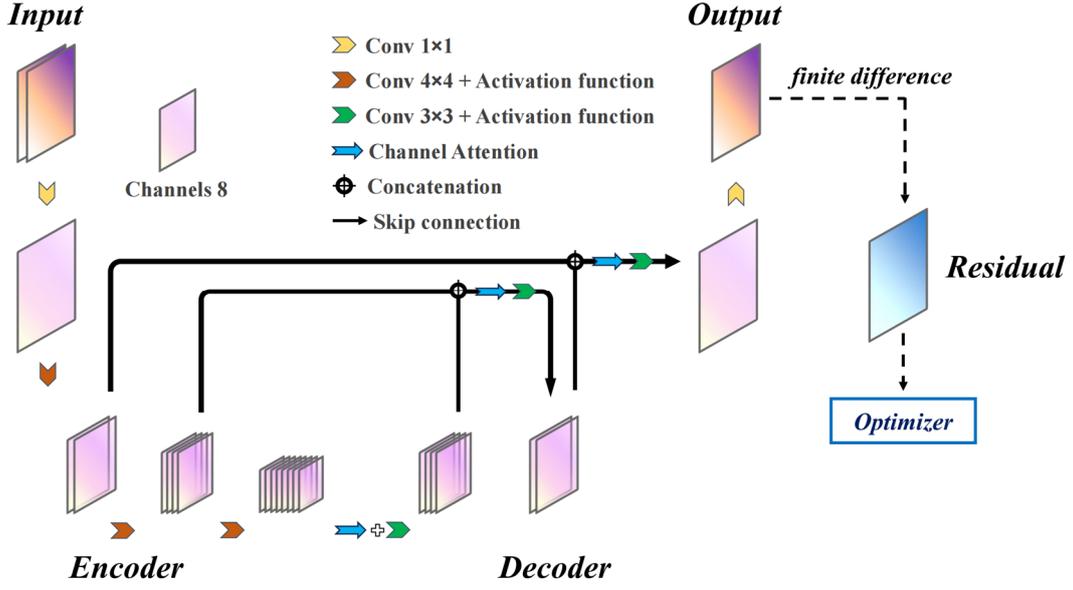

Figure 1. A schematic of PINNs, in which CNN is employed to establish a mapping from space coordinates to the solution.

2.2 Brief introduction of matrix preconditioning techniques

Consider a large-scale linear system

$$Au = b \tag{3}$$

where $A \in \mathbb{R}^{N \times N}$ and $b \in \mathbb{R}^{N}$. $A$ determines the properties of the system. Matrix preconditioning techniques aim to design an effective matrix $M \approx A$, known as the preconditioner, to improve the accuracy or reduce the computational cost of obtaining a numerical solution. A linear preconditioner is typically applied either on the left or the right, if suitably factored, as shown below.

$$\begin{aligned} M^{-1}Au &= M^{-1}b \\ AM^{-1}y &= b, u = M^{-1}y \end{aligned} \tag{4}$$

in Eq. (4), left preconditioning scales the residual vector, the condition number $\kappa(M^{-1}A)$ of $M^{-1}A$ is much smaller than that $\kappa(A)$ of $A$. Right preconditioning scales the solution vector. When substituted with $f = Au - b = 0$, left preconditioning can be expressed as $M^{-1}f = 0$. For more details on matrix preconditioning, please refer to [37, 38].

## 2.3 Matrix preconditioning for PINNs based on adjoint method

As stated in Section 2.1, the PDE is discretized using the finite difference to compute the loss function Eq. (2). As a result, $\mathcal{N}[\boldsymbol{u}^{NN}]=0$ is replaced by a discrete system $\boldsymbol{f}(\boldsymbol{u}^{NN})=\boldsymbol{0}$. For the linear PDE, $\boldsymbol{f}(\boldsymbol{u}^{NN})$ is a linear system:

$$\boldsymbol{f}^{NN} = \boldsymbol{A}\boldsymbol{u}^{NN} - \boldsymbol{b} = \boldsymbol{0} \tag{5}$$

where $\boldsymbol{A}$ is the Jacobian matrix $\boldsymbol{J}$:

$$\boldsymbol{J} = \frac{\partial \boldsymbol{f}^{NN}}{\partial \boldsymbol{u}^{NN}} \tag{6}$$

After appropriate linearization, the nonlinear PDE can also be represented as Eq. (5).

Thus, the loss function Eq. (2) becomes:

$$\mathcal{L} = \frac{1}{N_r}(\boldsymbol{f}^{NN})^T \boldsymbol{f}^{NN} \tag{7}$$

in our previous study [17], we verified the strong connection between the ill-conditioning of PINNs and the Jacobian matrix $\boldsymbol{J}$ by constructing a controlled system. Specifically, the training of PINNs becomes increasingly difficult or even fails as the condition number of the Jacobian matrix increases. Based on this conclusion, we propose a matrix preconditioning method for PINNs in this study, named Pre-PINNs. The left preconditioning is employed for the linear system obtained from the discretized PDE. Assuming the preconditioner is $\boldsymbol{M}$, the linear system after preconditioning becomes:

$$\boldsymbol{p}^{NN} = \boldsymbol{M}^{-1}\boldsymbol{f}^{NN} = \boldsymbol{M}^{-1}(\boldsymbol{J}\boldsymbol{u}^{NN} - \boldsymbol{b}) = \boldsymbol{0} \tag{8}$$

The Jacobian matrix of the linear system Eq. (8) is $\boldsymbol{M}^{-1}\boldsymbol{J}$. According to preconditioning theory, $\kappa(\boldsymbol{M}^{-1}\boldsymbol{J}) \ll \kappa(\boldsymbol{J})$. Thus, constructing the loss function based on Eq. (8) can improve the convergence and accuracy of PINNs. The new loss function is:

$$\mathcal{L} = \frac{1}{N_r}(\boldsymbol{p}^{NN})^T \boldsymbol{p}^{NN} \tag{9}$$

A key is obtaining the preconditioner $\boldsymbol{M}$, which is typically derived through factorization of the Jacobian matrix $\boldsymbol{J}$, such as the incomplete LU factorization used in this study. Thus, according to Eq. (6), we first compute $\boldsymbol{J}$ by combining automatic differentiation [13] with matrix coloring [39]. The purpose of introducing matrix coloring is to reduce memory and computational costs, which is crucial for large-scale problems.

Although $M$ is a sparse matrix, its inverse $M^{-1}$ is usually not sparse. For small-scale problems, $M^{-1}$ is available, which allows us to compute the loss using Eq. (9). However, for large-scale problems, obtaining $M^{-1}$ requires extremely high memory and computational costs. In practice, $p^{NN} = M^{-1} f^{NN}$ is often substituted by solving $Mp^{NN} = f^{NN}$ for $p^{NN}$ using iterative methods. For instance, with the incomplete LU factorization $J \approx LU = M$, where $L$ is lower unitriangular and $U$ is upper triangular. One then solves $Lq^{NN} = f^{NN}$ and $Up^{NN} = q^{NN}$, which can be done efficiently because the matrices are triangular. However, these solving processes are typically not supported by automatic differentiation. As a result, $p_{de}^{NN}$ is detached from the computational graph (subscript $de$ indicates the variable is detached from the computational graph), which makes it impossible to compute the gradients of the loss function with respect to the network parameters. Consequently, gradient descent cannot proceed. We introduce the adjoint method to address this issue.

We consider a constrained optimization problem:

$$\min_{f} \ (p_{de}^{NN})^T p_{de}^{NN}$$
$$\text{s.t.} \ \ p_{de}^{NN} = M^{-1} f_{de}^{NN} \tag{10}$$

To transform this constrained optimization problem into an unconstrained one, we construct the Lagrange function:

$$\mathcal{S} = (p_{de}^{NN})^T p_{de}^{NN} + \lambda^T (p_{de}^{NN} - M^{-1} f_{de}^{NN}) \tag{11}$$

where $\lambda$ is the Lagrange multiplier or the adjoint vector. The necessary optimality conditions for Eq. (11) are:

$$\frac{\partial \mathcal{S}}{\partial \lambda} = p_{de}^{NN} - M^{-1} f_{de}^{NN} = 0$$
$$\frac{\partial \mathcal{S}}{\partial p_{de}^{NN}} = 2 p_{de}^{NN} + \lambda = 0 \tag{12}$$
$$\frac{\partial \mathcal{S}}{\partial f_{de}^{NN}} = -(M^{-1})^T \lambda$$

where $p_{de}^{NN} - M^{-1} f_{de}^{NN} = 0$ is naturally satisfied because $p_{de}^{NN}$ is obtained by iteratively solving $Mp_{de}^{NN} = f_{de}^{NN}$. $\lambda = -2 p_{de}^{NN}$ is easily obtained from $2 p_{de}^{NN} + \lambda = 0$. Thus, we can compute the gradients $g_{de} = \partial \mathcal{S} / \partial f_{de}^{NN} = 2(M^{-1})^T p_{de}^{NN}$ of $(p_{de}^{NN})^T p_{de}^{NN}$ with respect to $f_{de}^{NN}$. Since $M^{-1}$ is not explicitly computed, we compute $g_{de}$ in the same way as $p_{de}^{NN}$, using the iterative method based on the incomplete LU factorization. Then, we construct a new loss function:

$$\mathcal{L} = \frac{1}{N_r}[(p_{de}^{NN})^T p_{de}^{NN} + g_{de}^T(f^{NN} - f_{de}^{NN})] \tag{13}$$

The value and the gradients with respect to the network parameters of Eq. (13) are the same as those in Eq. (9), and the former is more general. Thus, we use Eq. (13) as the loss function for Pre-PINNs.

## 3 Results

We validate the effectiveness of Pre-PINNs by solving Poisson equation and Navier-Stokes equations. To evaluate the accuracy of the solution, we compute the relative $L_2$ error:

$$\frac{\|u^{NN} - u^{Ref}\|}{\|u^{Ref}\|} \tag{14}$$

where $u^{Ref}$ represents the reference solution. It is analytical for Poisson equation and obtained using the finite volume method for Navier-Stokes equations. The study develops programs based on the PyTorch platform and executes the algorithm on NVIDIA GeForce RTX 3080 GPU. The code for this study is publicly available at https://github.com/Songstc/Pre-PICNNs.

3.1 Poisson equation

We first test the proposed method by solving the Poisson equation, which is a representative linear partial differential equation widely used in numerical experiments related to PINNs [14, 25, 40].

$$\frac{\partial^2 u}{\partial x^2} + \frac{\partial^2 u}{\partial y^2} = f(x,y), \quad (x,y) \in [-1,1]^2 \tag{15}$$

where the source term $f(x,y) = -\pi^2(k^2+1)\sin(k\pi x)\sin(\pi y)$, the exact solution $u(x,y) = \sin(k\pi x)\sin(\pi y)$. We set $k = 15$, according to the investigation by Basir and Senocak [41], PINNs fail in this case. Considering the high-frequency of $u$ in the $x$-direction, we use a uniform $N_x \times N_y = 300 \times 40$ mesh as the input to the CNN. Relying on the exact solution, the Dirichlet boundary conditions are enforced by padding. We solve this problem using PINNs and Pre-PINNs, respectively. The network input is $[x, y]$, the output is $u$, including three encoding layers and three decoding layers. Tanh is chosen as the activation function, and Adaptive Moment Estimation (Adam) optimizer is used to perform gradient descent. We train both networks for 50000 epochs. For Pre-PINNs, since the problem is linear, the Jacobian matrix $J$ only depends on the finite difference scheme and remains unchanged during the training process. Thus, before training Pre-PINNs, we compute it and further obtain the preconditioner $M$ through the incomplete LU factorization.

Due to the integration of the preconditioning method, the loss of Pre-PINNs is much smaller than that of PINNs. To enable a more reasonable comparison, we multiply the loss of Pre-PINNs by a weight of $w = 10^7$ in advance to ensure that the initial losses of both methods is on the same order of magnitude. The corresponding convergence history is shown in Figure 2. We observe that although the magnitudes of the loss reductions for both methods are comparable, the relative $L_2$ error of the results obtained by Pre-PINNs (1.38%) is more than an order of magnitude lower than that of PINNs (23.08%). Similar phenomena have also been observed in studies [7, 17, 42]. According to existing studies, due to issues such as the ill-conditioning of the PDE loss, the reduction in the PINNs loss function does not necessarily indicate more accurate results. It may fall into a non-physical local optimum. By reducing the condition number of the Jacobian matrix, Pre-PINNs is beneficial in avoiding this issue to some extent. Figure 3 shows the point-wise absolute error of the solutions obtained by both methods.

3.2 Navier-Stokes equations

Next, we validate the Pre-PINNs by solving incompressible Navier–Stokes equations:

$$\begin{aligned} &\frac{\partial u}{\partial x} + \frac{\partial v}{\partial y} = 0 \\ &u\frac{\partial u}{\partial x} + v\frac{\partial u}{\partial y} + \frac{\partial p}{\partial x} - \frac{1}{Re}(\frac{\partial^2 u}{\partial x^2} + \frac{\partial^2 u}{\partial y^2}) = 0 \\ &u\frac{\partial v}{\partial x} + v\frac{\partial v}{\partial y} + \frac{\partial p}{\partial y} - \frac{1}{Re}(\frac{\partial^2 v}{\partial x^2} + \frac{\partial^2 v}{\partial y^2}) = 0 \end{aligned} \qquad (16)$$

where $u$ is the velocity component in the $x$-direction, $v$ is the velocity component in the $y$-direction. $p$ is the pressure, $Re$ is the Reynolds number.

We study the lid-driven cavity flow on the computational domain $(x, y) \in [0,1]^2$, with boundary conditions as shown in Figure 4. Despite its simple geometry, the lid-driven cavity flow retains a rich fluid flow physics manifested by multiple counter rotating recirculating regions on the corners of the cavity as $Re$ increases [43]. We solve the lid-driven cavity flow for $Re = 500$ and $Re = 4000$ using a uniform $N_x \times N_y = 160 \times 160$ mesh, respectively. Since the flow at $Re = 4000$ is complex, we use the fourth-order central difference to compute the partial derivatives in Eq. (16). According to [17], PINNs fail to achieve any meaningful results when $Re = 4000$.

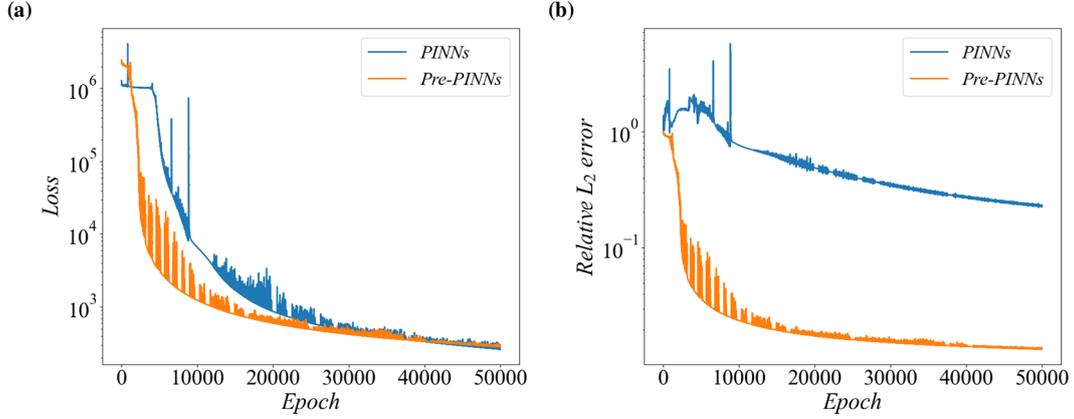

Figure 2. Comparison of the convergence history obtained by PINNs and Pre-PINNs with $w = 10^7$ for solving the poisson equation, respectively. (a) Loss. (b) Relative $L_2$ error.

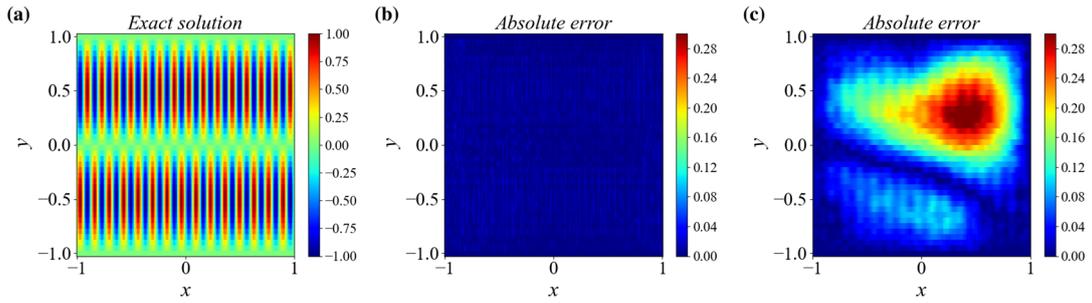

Figure 3. Point-wise absolute error between the exact solution and the results obtained by PINNs and Pre-PINNs with $w = 10^7$ for solving the poisson equation, respectively. (a) Exact solution. (b) Pre-PINNs. (c) PINNs.

We perform the linearization of the Navier-Stokes equations by transforming the convection terms $u\partial u/\partial x + v\partial u/\partial y$, $u\partial v/\partial x + v\partial v/\partial y$ into $u_{de}^{la}\partial u/\partial x + v_{de}^{la}\partial u/\partial y$, $u_{de}^{la}\partial v/\partial x + v_{de}^{la}\partial v/\partial y$, where $u_{de}^{la}, v_{de}^{la}$ are the network outputs from the previous epoch, and they are detached from the computational graph. We solve this problem using PINNs and Pre-PINNs, respectively. The network input is $[x, y]$, the output is $[u, v, p]$, including two encoding layers and two decoding layers. The Dirichlet boundary conditions are enforced by padding. GELU is chosen as the activation function. Limited-memory Broyden-Fletcher-Goldfarb-Shanno (LBFGS) optimizer is used to perform gradient descent, with the maximum number of iterations per epoch set to 200. We train PINNs for 1000 epochs, Pre-PINNs for 50 epochs under $Re = 500$, and for 100 epochs under $Re = 4000$. Same as Section 3.1, we multiply the loss function of Pre-PINNs by a weight to ensure that the initial losses of both methods are of the same order of magnitude. The weight $w$ is set to $10^3$ for $Re = 500$, and to $10^2$ for $Re = 4000$.

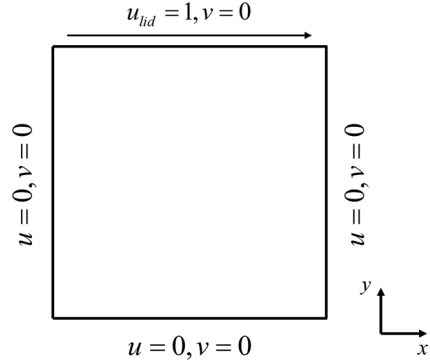

Figure 4. A schematic of boundary conditions for Lid-driven cavity flow.

Figure 5 shows the comparison of the convergence history obtained by PINNs and Pre-PINNs for solving Eq. (16) with $Re = 500$. We observe that the relative $L_2$ error of the result obtained by training Pre-PINNs for only 50 epochs ($0.93\%$) is nearly an order of magnitude lower than the result obtained by training PINNs for 1000 epochs ($8.12\%$), indicating that Pre-PINNs can improve both convergence speed and accuracy. The point-wise absolute error between the reference solution and the results obtained by PINNs and Pre-PINNs as illustrated in Figure 6. Figure 7 shows the comparison of the convergence history obtained by PINNs and Pre-PINNs for solving Eq. (16) with $Re = 4000$. We observe that although the loss function of the PINNs is gradually reducing, the relative $L_2$ error of its results is large all the time. In contrast, Pre-PINNs obtain results in good agreement with the reference solution after training for only 100 epochs, as shown in Figure 8. The potential of our method in alleviating the ill-conditioning of PINNs is demonstrated, which is crucial for applying PINNs to complex problems.

## 4 Conclusions

In this study, we propose a matrix preconditioning method to improve the convergence of PINNs, referred to as Pre-PINNs. Pre-PINNs modify the loss function through the preconditioner to reduce the condition number of the PDE system's Jacobian matrix, which is crucial for alleviating the ill-conditioning of PINNs. By combining matrix coloring and automatic differentiation, we first design an efficient method to compute the sparse Jacobian matrix. Subsequently, we employ the incomplete LU factorization for the Jacobian matrix to obtain the preconditioner and use the inverse of the preconditioner to scale the PDE residual in the loss function. Considering that computing the inverse of the preconditioner requires high memory and computational costs for large-scale problems, we use an iterative method to replace this process. Since the iterative method involving sparse matrices does not support automatic differentiation, we also develop a framework based on the adjoint method to compute the gradients of the modified loss function with respect to the network parameters. By solving the multi-scale problem governed by the Poisson equation and the high Reynolds number problem governed by the Navier–Stokes equations, in both of which PINNs fail, we demonstrate the effectiveness of Pre-PINNs.

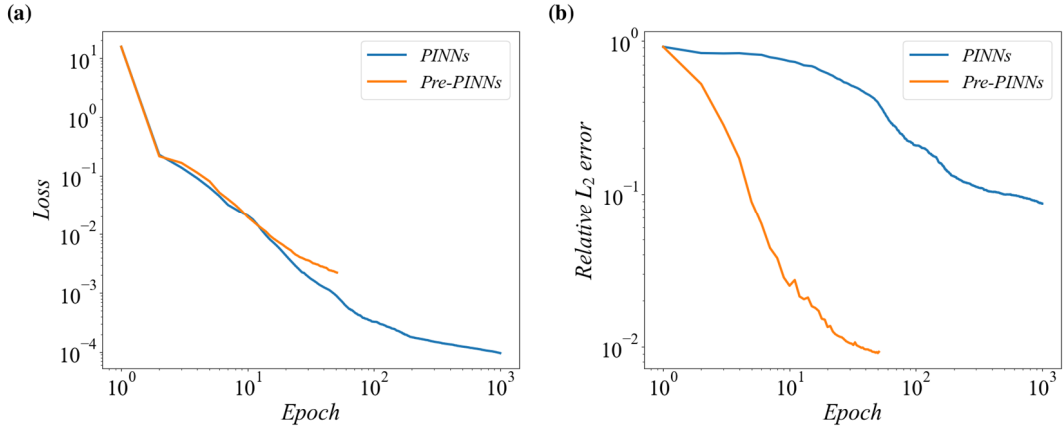

Figure 5. Comparison of the convergence history obtained by PINNs and Pre-PINNs with $w = 10^3$ for solving the lid-driven cavity flow with $Re = 500$, respectively. (a) Loss. (b) Relative $L_2$ error.

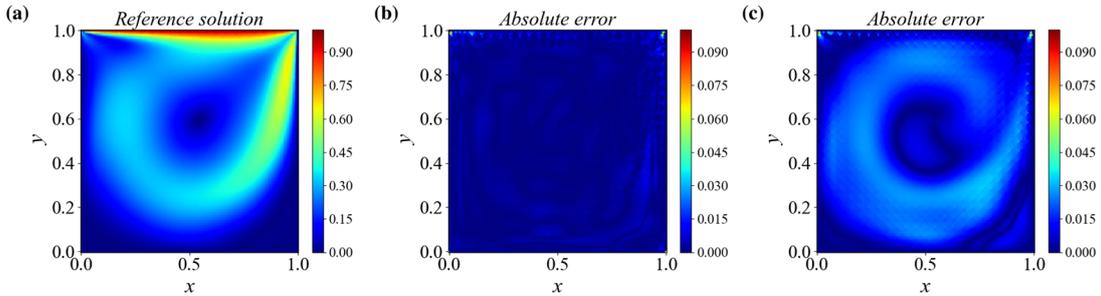

Figure 6. Point-wise absolute error between the reference solution and the results obtained by PINNs and Pre-PINNs with $w = 10^3$ for solving the lid-driven cavity flow with $Re = 500$, respectively. (a) Reference solution. (b) Pre-PINNs. (c) PINNs.

In the future, we will explore more applications of the proposed method. In addition, despite the advantages, Pre-PINNs still have some limitations. For instance, it relies on mesh for numerical discretization. We hope to develop a matrix preconditioning method suitable for meshfree approaches. One promising direction is to approximate the Jacobian matrix by leveraging the historical outputs of the neural network.

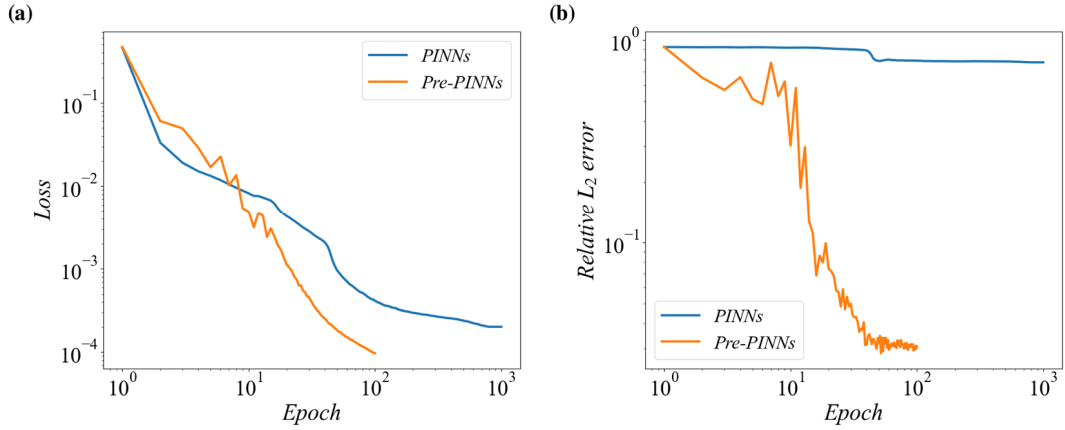

Figure 7. Comparison of the convergence history obtained by PINNs and Pre-PINNs with $w=10^2$ for solving the lid-driven cavity flow with $Re=4000$, respectively. (a) Loss. (b) Relative $L_2$ error.

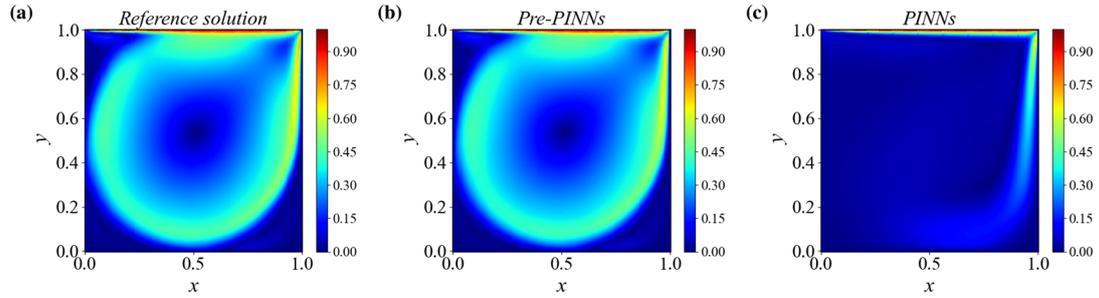

Figure 8. Comparison of the reference solution, the results obtained by PINNs and Pre-PINNs with $w=10^2$ for solving the lid-driven cavity flow with $Re=4000$, respectively. (a) Reference solution. (b) Pre-PINNs. (c) PINNs.

## Acknowledgements

We would like to acknowledge the support of the National Natural Science Foundation of China (No. 92152301).

## Conflict of Interest Statement

The authors have no conflicts to disclose.

## Data Availability Statement

The data that support the findings of this study are available from the corresponding author upon reasonable request.